\documentclass[12pt,a4paper]{article}
\usepackage[latin2]{inputenc}
\usepackage[T1]{fontenc}
\usepackage[leqno]{amsmath}
\usepackage{amsfonts,amssymb,latexsym,exscale}
\usepackage{graphicx}
\usepackage{pgf}
\usepackage{xspace}
\usepackage{theorem}
\usepackage{hyperref}
\usepackage{floatflt}
\usepackage[absolute]{textpos}
\usepackage{mathptmx}
\usepackage{textcomp}

\theoremstyle{plain}
\newtheorem{theorem}{Theorem}

\newtheorem{lemma}{Lemma}
\newtheorem{proposition}{Proposition}
\newtheorem{conjecture}{Conjecture}

\begin{document}

\title{\LARGE{\bf{ Optimal higher-dimensional Dehn functions\\ for some CAT(0) lattices}}}
\author{\Large{Enrico Leuzinger}}
\date{May 12, 2012}
\maketitle

\begin{abstract} Let  $X=S\times E \times B$ be the metric product of  a symmetric space $S$ of noncompact type,
 a Euclidean space $E$ and  a product $B$ of  Euclidean buildings.  Let $\Gamma$ be a discrete group acting isometrically and cocompactly on $X$. We determine a family of quasi-isometry invariants for such  $\Gamma$, namely the $k$-dimensional Dehn functions, which measure the difficulty to fill $k$-spheres by $(k+1)$-balls (for $1\leq k\leq \dim\ X-1$).   Since the group $\Gamma$ is quasi-isometric to the associated CAT(0)  space $X$,  assertions about Dehn functions for $\Gamma$ are equivalent to
 the corresponding results on filling functions for $X$.  Basic examples of groups $\Gamma$ as above  are uniform $S$-arithmetic subgroups of reductive groups defined over global fields. We also discuss a (mostly) conjectural picture for
non-uniform $S$-arithmetic groups.
\end{abstract}

\vspace{6ex}



\maketitle


\section{Introduction and main results}

Since Gromov's seminal essay \cite{Gr2} the study of quasi-isometry invariants of 
(Cayley graphs of) finitely generated groups is of central interest in geometric group theory.
 Dehn functions are basic and well-studied examples of such quasi-isometry invariants.
They measure  the 
difficulty to fill circles by discs and are closely related to the complexity of the word problem
(see e.g. \cite{Bri1}).
For example, the Dehn function of a CAT(0) group is either linear or quadratic. 
A finitely generated  group is Gromov-hyperbolic if and only if its Dehn function is linear (see \cite{Bri1}, 6.1.5). 

More generally,  higher-dimensional homotopical resp. homological Dehn (or filling) functions are quasi-isometry invariants which encode the difficulty to fill
 $k$-spheres by $(k+1)$-balls resp. $k$-cycles by $(k+1)$-chains. They were first considered in \cite{E+} and \cite{Gr2} \S 5, and then more systematically  in \cite{AWP}. In contrast to the classical 1-dimensional Dehn functions
much less is known about higher-dimensional Dehn functions. Recently, versions of such quasi-isometry invariants have been investgated by several authors for various classes of finitely generated groups (see e.g. \cite{ABDDY},  \cite{BBFS}, \cite{Gft}, \cite{Ril}, \cite{Yo1}, \cite{Yo2}, \cite{Yo3}).

The chief goal of the present  paper is Theorem 2, in which we precisely determine the higher-dimensional Dehn functions for groups  $\Gamma$ acting geometrically 
(i.e., isometrically, properly discontinuously and cocompactly) on a
product $X=S\times E \times B$, where $S=G/K$ is a symmetric space of noncompact type
(with $G$ a semisimple, noncompact real Lie group and $K\subset G$ a maximal compact subgroup),
$E$ is a Euclidean space and $B$ is a product of Euclidean buildings. 
These higher-dimensional Dehn functions   are (essentially) equivalent to the corresponding filling functions of the associated complete simply connected CAT(0) space $X$, which we compute in Theorem 1.

 Examples of groups  $\Gamma$ for which Theorem 2 holds are crystallographic groups, uniform lattices
(arithmetic or not)  in semisimple real Lie groups and (uniform) $S$-arithmetic subgroups of 
connected, reductive and $\mathbb K$-{\it anisotropic} $\mathbb K$-groups (where $\mathbb K$ is a global field). We briefly discuss this last family in Section 4. Our results are motivated by the corresponding problem to determine the higher-dimensional Dehn functions for $S$-arithmetic subgroups of 
connected, reductive and $\mathbb K$-{\it isotropic} $\mathbb K$-groups, which yield \emph{non-uniform} lattices acting on CAT(0) spaces  $X$ as above. Since such a non-uniform lattice is not quasi-isometric to $X$, this
 latter problem is much more difficult. In Section 5 we solve it  for arbitrary lattices in $SO(n,1)$ and $SU(n,1)$ (see Theorems \ref{hyperbol} and 5). We also  provide  a (mostly) conjectural picture for the general case (see Conjecture 1).

Beside  Heisenberg groups studied by Young in \cite{Yo2}, the groups $\Gamma$ covered by Theorems 2, 4 and 5 seem currently to be
 the only examples of families  of finitely generated discrete groups  whose Dehn functions are  known \emph{in all} dimensions.

\subsection{Higher-dimensional filling functions for some CAT(0) spaces}

In order to state our results we first recall the concept of higher-dimensional filling functions for
CAT(0) spaces  $X=S\times E \times B$ as above. Consider an (integral) Lipschitz $k$-chain in $X$, i.e., a finite linear combination
$\Sigma=\sum_i a_i\sigma_i$, with $a_i\in \mathbb Z$ and $\sigma_i:\Delta^k\to X$  a Lipschitz map
from the $k$-dimensional Euclidean  standard simplex  to $X$. Note that $S$ and $E$ are smooth manifolds and $B$ is a piecewise Euclidean (poly)-simplicial complex with bounded geometry. Therefore, by Rademacher's theorem
(\cite{Fed}, 3.1.6), such a $\sigma_i$ is differentiable 
almost everywhere and  we thus can define $\textup{vol}_k\sigma_i$ as the integral over its Jacobian
(for smooth maps this is the Riemannian volume of the image, compare \cite{BBI}, 5.5). We then define the $k$-\emph{volume} or $k$-\emph{mass} of $\Sigma$ as 
$$\textup{vol}_k\, \Sigma:=\sum_i|a_i|\textup{vol}_k\sigma_i.$$

 We wish to  measure the difficulty to fill Lipschitz $k$-cycles by Lipschitz $(k+1)$-chains. More precisely,  for  an integral Lipschitz $k$-cycle $\Sigma$ we define its \emph{filling volume}
$$
\textup{FVol}_{k+1}(\Sigma):=\inf\{\textup{vol}_{k+1}\, \Omega\mid \Omega=\textup{Lipschitz}\  (k+1)\textup{-chain with}\ \partial \Omega=\Sigma\}.
$$
Notice that the CAT(0) space $X=S\times E \times B$ is contractible  and thus  $k$-connected for any $0\leq k\leq n-1=\dim X-1$ (i.e., all  homotopy groups $\pi_k(X)$ are trivial); hence  the corresponding homology groups are also trivial (by the Hurewicz isomorphism theorem) and such fillings exist.

The \emph{ $(k+1)$-dimensional 
filling function} of $X$ is  then  given by
$$
\textup{FV}_X^{k+1}(l):=\sup\{\textup{FVol}_{k+1}(\Sigma)\mid \Sigma=\textup{Lipschitz}\  k\textup{-cycle in X with}\ \textup{vol}_k\, \Sigma\leq l\}.
$$

We emphasize that we are interested in the asymptotic geometry of $X$, i.e., the cycles we wish to fill
are supposed to be ``large''. In fact, volumes of ``small'' domains have
 Euclidean behaviour and in particular satisfy superlinear isoperimetric inequalities
(see \cite{Gr3}, p. 322).

Quasi-isometry invariants are obtained from  filling functions by considering equivalence classes for
the following (standard) equivalence relation.
 Let $f,g:\mathbb R\to \mathbb R$ be two functions. We  write $f\prec g$,
if there is a constant $C>0$ such that $f(x)\leq Cg(Cx+C)+Cx+C$ and $f\sim g$,
if $f\prec g$ and $g\prec f$. 
Our key result is the following 

\begin{theorem}  \textup{(Filling functions)}\ \ Let $X=S\times E \times B$ be the metric product of  a symmetric space of noncompact type $S=G/K$,
a Euclidean space $E$ and a product  $B$ of Euclidean buildings (one or two factors of $X$ may be trivial). Let $\textup{rank}\ X=\textup{rank}_{\mathbb R} G+\dim E+\dim B$ be the Euclidean rank of $X$. For $k\geq 1$, let $\textup{FV}_X^{k+1}$ be the $(k+1)$-dimensional filling volume function of $X$.
Then

\noindent \textup{(i)} \ \  $X$ has  Euclidean  filling functions  below the rank:
$$\textup{FV}_X^{k+1}(l)\sim  l^\frac{k+1}{k}\ \ \textup{ if}\ \  k \leq \textup{rank}\ X -1;$$
\textup{(ii)} \ \  $X$ has linear filling functions above the rank:
$$\textup{FV}_X^{k+1}(l)\sim l\ \ \textup{ if}\ \  \textup{rank}\ X\leq k\leq \dim X-1.$$
\end{theorem}

 Notice that case (ii) in Theorem 1  only occurs if the symmetric factor $S$ is nontrivial.

If $X=S$ is a symmetric space of noncompact type, Theorem 1 has been correctly asserted by Gromov in \cite{Gr2},
$ 5.D(5)(b')$. He also proposed a 
possible proof for an \emph{upper bound} via isoperimetric inequalities: Pick some maximal flat $F$ in $X$ and  project othogonally to that flat. Then use 
the resulting cylinder to produce the desired filling. Gromov claims that this
projection exponentially retracts the volume (as a function of the distance $d$ to $F$). 
In $6.B_2 (d'')$ he attributes this to Mostow. But Mostow actually only proves (using quite 
delicate estimates)  that  the contraction
is proportional to ${d}^{-1/2}$ when restricted to an r-dimensional submanifold (see \cite{Mo}, Lemma 6.4).
The problem is that  vectors tangent  to singular geodesics can have projections with contraction factor almost equal to one, so that it is not clear whether this idea works. Nevertheless our approach is still quite
 similar: instead of projecting to a flat we project to a suitable horosphere.

\vspace{1ex}

A result of Wenger  \cite{We1} asserts that a complete, simply connected CAT(0) space satisfies  Euclidean isoperimetric inequalities
in all dimensions (see  Proposition 1 below). This yields Euclidean \emph{upper bounds} for filling functions  (see Section 2.1). Theorem 1 shows that this inequalities are sharp below the rank. Moreover it refines Wenger's general result in the present setting for dimensions above the rank (again showing that linear isoperimetric inequalities are sharp by providing lower bounds). In particular, Theorem 1
asserts that the filling functions of a symmetric space  detect its rank.

Hyperbolic spaces satisfy linear isoperimetric inequalities in all dimensions.
 Theorem 1 thus also goes with  Gromov's philosophy concerning the notion of ``rank'' in non-positive curvature as expressed in \cite{Gr2}, $6.B_2$:
\begin{quote}One of the guiding principles in the asymptotic geometry of spaces $X$ with $K\leq 0$
 can be expressed as follows. All flatness of $X$ where $K=0$ is confined to $k$-flats in $X$ ... One 
distinguishes among them \emph{maximal} flats ... and then tries to show that $X$ is ``hyperbolic 
transversally to maximal flats''.
\end{quote}

\subsection{Higher-dimensional Dehn functions for some CAT(0) groups}

We now pass from spaces to groups. Recall that a group $\Gamma$ is of type $F_m$, if there exists an
Eilenberg-Mac Lane complex $K(\Gamma,1)$ with finite $m$-skeleton. A group is of type $F_{\infty}$, if it is of type $F_m$ for all $m$. The finiteness properties 
$F_1$ and $F_2$ are equivalent to $\Gamma$ beeing finitely generated and finitely presented, respectively.
Alonso showed in [Al] that ``being of type $F_m$'' is  invariant under quasi-isometries.

Suppose that a given group $\Gamma$ is of type $F_{k+1}$. Then there exists a
 $k$-connected $(k+1)$-complex $K$ on which $\Gamma$ acts cocompactly (see \cite{Bri2} and
\cite{Ril}). A combinatorial map  $\Sigma: S^k\to K^{(k)}$ of the standard $k$-sphere  into the $k$-skeleton of $K$ can then be filled by a combinatorial $(k+1)$-ball $\Omega$. 
We define the combinatoral  $k$-volume of $\Sigma$ as the number of non-degenerate $k$-cells in the image
$$
\textup{vol}^c_{k}(\Sigma):=\sharp \{\sigma\subset \Sigma(S^k)\mid \sigma \, \textup{is a non-degenerate $k$-cell of}\  K\}.
$$
Analogously we define the combinatorial $(k+1)$-volume of $\Omega$ as the number of non-degenerate $(k+1)$-cells in the image of $\Omega$.
We further  set
$$
\textup{DFVol}^c_{k+1}(\Sigma):=\min\{\textup{vol}^c_{k+1}(\Omega)\mid \Omega\ \textup{is an $(k+1)$-ball with boundary}\ \Sigma\}.
$$

The \emph{$k$-dimensional (combinatorial) Dehn function} (of $K$) is then  defined as
$$
\delta^{k}_{K}(l):= \max\{
\textup{DFVol}^c_{k+1}(\Sigma)\mid   \textup{vol}^c_{k}(\Sigma)\leq l\}.
$$

The above definition of higher dimensional Dehn functions is equivalent to  that of Alonso et al. in \cite{AWP} (see \cite{Bri2}).
In particular its equivalence class (in the sense above) is independent of the chosen complex $K$ and  defines a quasi-isometry invariant of $\Gamma$.

  \vspace{1ex}

We next give an alternative definition of higher-dimensional Dehn functions which is better adapted to 
our situation, since it is more closely related to filling functions.

Let $X$ be a $k$-connected manifold or Euclidean polyhedral  complex and let  $\Gamma$ be a (finitely generated) group which  acts cocompactly on $X$.
The \emph{$k$-th dimensional (Lipschitz) Dehn function} (of $X$), $\delta_X^k$, is a homotopical version of a filling function. It measures the volume that is  necessary to extend a Lipschitz map $f: S^k\to X$ of the $k$-dimensional unit sphere 
to a map $g: B^{k+1}\to X$ of the $(k+1)$-dimensional unit ball:
$$
\textup{DFVol}_{k+1}(f):= \inf\{\textup{vol}_{k+1}(g)\mid g:B^{k+1}\to X,\  g_{\mid_{S^k}}=f\}
$$
and 
$$
\delta^k_X(l):= \sup\{
\textup{DFVol}_{k+1}(f)\mid f:S^k\to X,\  \textup{vol}_k(f)\leq l\}.
$$

 As a consequence of the Federer-Fleming deformation theorem the combinatorial and Lipschitz versions of higher dimensional Dehn functions are actually  equivalent if $X$ (or $K$) 
is a space which can be approximated by a polyhedral complex of bounded geometry  (compare \cite{ABDDY}, Theorem 2.1.). 
Since its equivalence class defines a  quasi-isometry invariant (by \cite{AWP}), we will
  write $\delta^k_{\Gamma}$ in place of $\delta^k_X$ or $\delta^k_K$  to denote (a representative of) this equivalence class.

\vspace{1ex}

 Theorem 1 implies

\begin{theorem}  \textup{(Dehn functions)}\ \ Let $X=S\times E \times B$ be the metric product of  a symmetric space of noncompact type $S=G/K$,
a Euclidean space $E$ and a product  $B$ of Euclidean buildings (one or two factors of $X$ may be trivial) and let $\textup{rank}\ X=\textup{rank}_{\mathbb R} G+\dim E+\dim B$ denote the Euclidean rank of $X$.

 Let $\Gamma$ be a group acting  isometrically,  cocompactly and properly discontinuously on X. 
Then

\textup{(i)}   the $k$-dimensional Dehn functions of $\Gamma$ are Euclidean below the rank:
$$
\delta_{\Gamma}^{k}(l)\sim l^{\frac{k+1}{k}}\ \ \ \textup{ if} \ \  k\leq \textup{rank}\ X-1\ ;
$$

\textup{(ii)}  the $k$-dimensional Dehn functions of $\Gamma$ are linear above the rank:
$$
\delta_{\Gamma}^{k}(l)\sim l\ \ \ \textup{ if} \ \ \ \textup{rank}\ X\leq k\leq \dim X -1.
$$ 
\end{theorem}

As noted above, the chief examples of groups $\Gamma$ to which Theorem 2 applies are provided by cocompact S-arithmetic subgroups of reductive algebraic groups $\mathbf G$ which are defined over  global fields $\mathbb K$. We will discuss such groups in some more detail in Section 4. Notice that, if  $\mathbb K$ is a function field, then Theorem 2 reduces to part (i) only.
Note further that in the special case where  $X=S$ is a symmetric space of \emph{rank $1$},
a group $\Gamma$ as in Theorem 2 is hyperbolic and then Theorem 2 asserts that all Dehn functions are linear.
For general hyperbolic groups linear upper bounds were obtained in  \cite{La} and \cite{Mi}. 

As already mentioned, linear $1$-dimensional Dehn functions characterize Gromov-hyper\-bolic groups.
 On the other end of the range, top-dimensional linear isoperimetric inequalities (or, equivalently, 
 positive Cheeger constants) characterize discrete \emph{non-amenable} groups, 
which act geometrically (see \cite{Gr2}, 05.A6 and Theorem \ref{topdim} below). 

 Finally, we note that a theorem of Avez \cite{Av} asserts that
the fundamental group of a closed manifold of nonpositive sectional curvature $V$ is non-amenable unless $V$ is flat.
This in particular applies to fundamental groups of compact locally symmetric spaces
 (of noncompact type)  with Euclidean factors, i.e., uniform lattices in real reductive Lie groups.
Theorems 1 and 2  thus refine  Avez's theorem in the case of reductive locally symmetric spaces.

\vspace{2ex}

{\bf Remark}\  (Generalization)\ Theorem 1 and Theorem 2  extend to CAT(0) spaces (and groups acting geometrically on them) of the form $X=S\times H$, where $S$ is a symmetric space of noncompact type as before, and where  $H$ is a Euclidean polyhedral complex  with
bounded geometry and  the property that its dimension is equal to its Euclidean rank  (i.e., the maximal dimension of an isometrically embedded Euclidean space).  In fact, in the proofs below we  only use
these (weak) geometric properties of buildings.

\section{ The proof of Theorem 1} 

\subsection{Theorem 1(i): Euclidean fillings below the rank}

By definition of the filling functions $\textup{FV}_X^{k+1}$ we have $k\geq 1$. As we are here concerned with Theorem 1(i), we may  assume  that  $\textup{(Euclidean) rank}\ X\geq 2$.

Wenger proved Euclidean isoperimetric inequalities in a general setting:

\begin{proposition} \ \textup{(Wenger, \cite{We1})}\ \label{cat0-iso} If $X$ is a complete simply connected CAT(0) space  of dimension $n$, then any $k$-dimensional Lipschitz cycle $\Sigma$ in $X$  satisfies a  Euclidean isoperimetric inequality:
$$
\textup{FVol}_{k+1}(\Sigma)\leq C_k\ (\textup{vol}_k\, \Sigma)^{\frac{k+1}{k}},\ \ 1\leq k\leq n-1,
$$
where $C_k$ is a constant that depends only on $X$ and $k$.
\end{proposition}

 These isoperimetric inequalities provide an upper bound for filling functions
$$\textup{FV}_X^{k+1}(l)\prec l^{\frac{k+1}{k}}.$$

 In order to get also a lower bound 
for $k<r=\textup{rank}\, X$ we consider a $k$-dimensional round sphere $S$ of volume $l$ contained  in a maximal ($r$-dimensional) flat $F$ of $X$.
Let $B$ be a minimal filling of $S$ \underline{in $X$}. The orthogonal projection from $X$ to $F$ is a $1$-Lipschitz map
(see e.g. \cite{BH}, 2.4).
If we  orthogonally project $B$ to $F$ we thus get a 
filling $B'$ of $S$ \underline{in  $F$} with $\textup{vol}_{k+1}\, B'\leq \textup{vol}_{k+1}\, B$ (see \cite{BBI}, 5.2.2).
Hence $B'$ is a minimal filling of the round sphere $S$ in the the Euclidean space $F$. By the solution of the isoperimetric problem in Euclidean space, $B'$ must be  a round ball with boundary $S$ and thus
satifies $l^{\frac{k+1}{k}}\sim  \textup{vol}_{k+1}\, B' \leq \textup{vol}_{k+1}\, B$.
Hence we get the required lower bound $$l^{\frac{k+1}{k}}\prec \textup{FV}_X^{k+1}(l).$$

\subsection{Theorem 1(ii): Linear fillings above the rank}

Note that Theorem 1 (ii), which concerns  fillings of dimensions above rank $X$, is  only meaningful if the symmetric factor $S$ of $X$ is nontrivial. In fact, otherwise we have  $X=E\times B$ and
the (Euclidean) rank of $X$ coincides with the topological dimension $\dim X= \dim E+\dim B$ and hence there
is no filling problem above the rank. \emph{In this subsection we will thus always assume that the symmetric factor $S$ is nontrivial}.

\subsubsection{Linear upper bounds}

We first recall the formulae for Jacobi fields in  Riemannian
products   $S\times M$, where $S=G/K$ is a symmetric space of noncompact type (with  semisimple Lie group $G$) and $M$ is an open (flat) submanifold of some Euclidean space. These will be used below for volume estimates. A proof for semisimple symmetric spaces can be found in \cite{Eb}, 2.14.
and 2.15; it directly extends to the product case. Recall that a Jacobi field along a geodesic ray is called 
\emph{stable}, if its norm is
bounded.

\begin{proposition}\label{jac-sym} Let $S\times M=G/K\times  M$ be an $n$-dimensional Riemannian product of a symmetric space of noncompact type with an open (flat) submanifold $M\subset \mathbb R^m$ of a Euclidean space and let  $x_0$ be a base point of $S$. Let $\mathfrak{g}=\mathfrak{k}\oplus \mathfrak{p}$
 be a Cartan decomposition of the Lie algebra of the semisimple group $G$. Let $\mathfrak{a}$ be a maximal abelian
 subspace of $\mathfrak{p}$. For a unit vector $ H$ in
$\mathfrak{a}$, an element $g\in G$ and a point $p\in M$, let $$c:[0,\infty)\to S\times M, \ \ t\mapsto (g\cdot\exp tH\cdot x_0, p)$$ be a (singular) geodesic ray parametrized by arc-length starting at $(g\cdot x_0,p)$. Further let
 $E_l(t), 0\leq l\leq n-1 $, be orthonormal parallel vector fields in $S\times M$ along the ray $c$ and
 orthogonal to $c$. 
Then any stable Jacobi field $Y(t)$ along $c(t)$ and orthogonal to $c$ with $Y(0)=\sum_{l=1}^{n-1}y_lE_l(0)$ is of the form
$$
Y(t)=\sum_{l=1}^{n-1}y_le^{-\sqrt{\lambda_l}\, t}E_l(t),
$$
where the $\lambda_l$ are zero or eigenvalues of the curvature operator 
$R_{H}=(\textup{ad}H)^2_{\mid \mathfrak{p}}$ and thus either also equal
 to zero or of the form $\lambda_l=(\alpha_l(H))^2$
for some (positive) root $\alpha_l$.
\end{proposition}

We now proceed to prove the linear upper bound in Theorem 1(ii). To that end let $X=S\times E\times B$
be as in Theorem 1 (with  nontrival symmetric factor $S$). Further let $\Sigma$ be a $k$-dimensional cycle with $$k\geq r=\textup{rank}\ X=\textup{rank}\ S+\dim\ M=\dim \ \mathfrak{a}+\dim\ M.$$ For convenience we use in the following  the same notation for the map $\Sigma$ and its image in $X$. We  wish to construct a chain 
$\Omega$
which bounds $\Sigma$ and such that $\textup{vol}_{k+1}\Omega\leq \textup{const}\cdot
\textup{vol}_k\Sigma$, where the constant only depends on $X$ and $k$.
The idea is to fill $\Sigma$ by first projecting it  orthogonally to a suitable horosphere of $X$ and 
 then take
$\Omega$ as  the cylinder of this map  together with a  filling of the projection 
(which will be very small, i.e., of volume $\leq 1$) inside this horosphere. We will achieve  this in 6 steps.

\vspace{1ex}

{\sc STEP 1 (A cone over $\Sigma$)}:\ Consider a Weyl chamber in the geometric boundary of 
the symmetric factor $S=G/K$ and denote its barycenter by $\omega$.
Together with the base point $x_0$ of $S$ this
 defines an Iwasawa decomposition $G=NAK$ and
 associated horocyclic coordinates $S\cong NA\cdot x_0$ (see e.g. \cite{Le}). 
Let $c_1(t):=\exp tH_0\cdot x_0,\  t\in \mathbb R,$ be the unit speed geodesic in $S$ through $x_0$ and asymptotic to $\omega$. Note that taking translates of $c_1$ gives a foliation 
of $S$ by asymptotic geodesics: $S=N\cdot \exp H_0^{\perp}\cdot c_1(\mathbb R)$. 

Given a point $x=(n\exp H_1 c_1(s),e,b)\in X=S\times E\times B$, with $H_1\in H_0^{\perp}$, we denote by $\tau(x,t)$ the unique (singular) geodesic ray in $X$  which starts at $x$ and is
asymptotic to the above barycenter $\omega\in S(\infty)\subset X(\infty)$, i.e.,
$$
\tau(x,t)=(n\exp H_1\cdot c_1(s+t),e,b), \ \  t\in [0,\infty).
$$
 This in particular defines  a map
$$
\tau: \Sigma\times \mathbb R_{\geq 0}\to X,\ \  (x,t)\mapsto \tau(x,t).
$$
Below we will use the  associated $(k+1)$-dimensional ``cone'' $\tilde{\Omega}:=
\tau(\Sigma\times \mathbb R_{\geq 0})\subset X$ over the $k$-cylce $\Sigma$.

In order to guarantee that $\tilde{\Omega}$ is $(k+1)$-dimensional, we pick some point $\tilde{x}\in \Sigma$ and choose the barycenter
$\omega$ in such such a way that the geodesic ray $\tau(\tilde{x},t)$ is \emph{transversal} to $\Sigma$.
That this is possible for any $k$ follows e.g. from the polar decomposition of $S$ (and its boundary at infinity):
\cite{Hel} Ch. 5, Lemma 6.3. We also note that the image under $\tau$ of the set of those points of 
$\Sigma$ where the issuing geodesic ray is tangential (i.e. not transversal) is only a $k$-dimensional subset and hence does not contribute to the  $(k+1)$-volume of the cone $\tilde{\Omega}$.

\vspace{1ex}

{\sc STEP 2 (Regular  points of $\Sigma$)}:\ 
If the building factor $B$ of $X$ is nontrivial, $X$ is not a differentiable manifold. We can, however, consider all open cells of $B$, say $\sigma_m$ with $m$ in some index set $I_B$. Then there is a finite subset $I\subset I_B$ such that the intersections 
$$\Sigma_m:=\Sigma\cap (S\times E\times \sigma_m),\  m\in I,$$ 
have positive $k$-dimensional Hausdorff measure and cover $\Sigma$ up to a set of measure zero (recall that we use here the same notation for the map $\Sigma$ and its image in $X$). 
 Notice that $S\times E\times \sigma_m, m\in I_B,$ 
is the Riemannian product of a symmetric space with an (open, flat) Euclidean manifold, so that Proposition 2 can be applied.

We fix an $\epsilon >0$, which will be chosen appropriately in STEP 5. 
By Rademacher's theorem  $\Sigma_m$ is differentiable almost everywhere
(with respect to Hausdorff measure). By Lusin's theorem (see \cite{Fed}, 2.3.5) there is  a compact subset $\Sigma_m^r$ of $\Sigma_m$ for each $m\in I$, such that $\Sigma_m^{r}$ is a $C^1$-manifold and $\textup{vol}_k(\Sigma_m\setminus
\Sigma_m^{r})<\epsilon$.

The sets $\Sigma_m^{r}$
are non-empty, disjoint $C^1$-manifolds which cover the regular part $$\Sigma^{r}:=\bigcup_{m\in I}\Sigma_m^{r}$$ of $\Sigma$. Note that  $\textup{vol}_k(\Sigma\setminus
\Sigma^{r})< \epsilon \mid\!I\!\mid$.

\vspace{1ex}

{\sc STEP 3 (Volume of the cone over the regular part $\Sigma^{r}$)}:\ 
We fix $m\in I$ and, for $\tau $ as ins STEP 1, we
denote by  $\omega$ be the volume form on $\tilde{\Omega}_m^{r}:=\tau(\Sigma_m^r\times \mathbb R_{\geq 0})\subset X$.
  We also define $\phi:\Sigma_m^{r}
\times\mathbb R_{\geq 0}$
by 
\begin{equation}\label{dens}
\tau^*\omega(x,t)=\phi(x,t)dt\wedge
\pi^*\bar{\omega}(x)
\end{equation}
where $\pi:\Sigma_m^r\times\mathbb R_{\geq 0}\to \Sigma_m^r$ is the projection to the first factor and 
$\bar{\omega}$ is the volume form on $\Sigma_m^r$.
 We will show that there is a constant  $C_2$ depending only on $X$ such that
\begin{equation}\label{localest1}
\textup{vol}_{k+1}(\tilde{\Omega}_m^{r})\leq C_2 \ \textup{vol}_{k}(\Sigma_m^r)\ \ \ \textup{for all}\  m\in I.
 \end{equation}
Equation (\ref{localest1})
 in turn implies the following linear estimate for the cone over the regular part $\Sigma^{r}$ of $\Sigma$
\begin{equation}\label{localest}
\textup{vol}_{k+1}(\tilde{\Omega}^{r})=\sum_{m\in I}\textup{vol}_{k+1}(\tilde{\Omega}_m^{r})
\leq C_2 \  \sum_{m\in I}\textup{vol}_{k}(\Sigma_m^r)=C_2\ \textup{vol}_{k}(\Sigma^{r}).
 \end{equation}

In order to prove estimate (\ref{localest1}) we will use the following

\begin{lemma}  Let $k\geq \textup{rank}\ X$,  $m\in I$ and $\epsilon$ as above. Then there are  constants
$\lambda_*>0$,  $C_1>0$ and $C_2>0$  depending only on $X$  (but not on $k$, $m$ or $\epsilon$) such that for
every $x\in \Sigma_m^r$ the density function $\phi$ in (\ref{dens}) satisfies
$$
\phi(x,t) \leq C_1e^{-\sqrt{\lambda_*}\, t}\ \  and\  thus \ \  \int_0^{\infty}\phi(x,t)dt= C_2<\infty .$$
 \end{lemma}

\emph{Proof}.\   Consider the geodesic ray $c(t):=\tau(x,t)$ defined in STEP 1 starting at $x\in \Sigma_m^r$ and asymptotic to $\omega$. We can (and do) assume that $c$ is transversal to $\Sigma$ (see the discussion in STEP 1). Let $V_0(0),\ldots,V_{n-1}(0)$
be an orthonormal frame of $T_xX$ such that $V_0(0)=\dot{c}(0)$ and $V_0(0),\ldots, V_{k}(0)$
span $T_x\tilde{\Omega}_m^{r}$. For each $0\leq i\leq k$ the unit vector $V_i(0)$ extends to
a parallel vector field $V_i(t)$ along $c(t)$. We choose local coordinates $(x_1,\ldots, x_{k})$
in $\Sigma_m^r$ around $x$ such that $\frac{\partial}{\partial x_i}(x)=V_i(0)+a_iV_0(0)$.
Then, by the map $\tau$, we have local coordinates $(x_1,\ldots, x_{k},t)$  in $\tilde{\Omega}_m^{r}$
 near the geodesic  ray $c(t)$. The volume form $\bar{\omega}$
of $\Sigma_m^r$ at $x$ is 
$$\bar{\omega}(x)=\sqrt{1+\sum_ia_i^2}\ dx_1\wedge\ldots\wedge dx_{k}.$$

We now consider the (unique) stable Jacobi fields $Y_i(t)$ along $c(t)$ which satisfy
 $Y_i(0)=V_i(0), 1\leq i\leq k$. Using the $k\times k$ matrix
$$
A_k(t):=(\langle Y_i(t),Y_j(t)\rangle)_{1\leq i,j\leq k}
$$
we can write the volume element $\omega$ on $\tilde{\Omega}_m^{r}$ as
$$
\omega(c(t))=\sqrt{\det A_k(t)}\,dt\wedge dx_1\wedge\ldots\wedge dx_k.
$$
Hence we get 
\begin{equation}\label{phi}
\phi(x,t)=\sqrt{\frac{\det A_k(t)}{1+\sum_ia_i^2}}\leq \sqrt{\det A_k(t)}.
\end{equation}

We now choose 
orthonormal parallel vector fields $E_l(t), 0\leq l\leq n-1,$ along the ray $c(t)$ as in Proposition \ref{jac-sym}
such that, for $ 1\leq i\leq n-1$,
$$
V_i(0)=\sum_{l=1}^{n-1}a_{il}E_l(0).
$$
Note that the matrix $(a_{il})$ is orthogonal.
By Proposition \ref{jac-sym}, the Jacobi fields $Y_i,  1\leq i\leq k,$ can  be written as
$$
Y_i(t)=\sum_{l=1}^{n-1}e^{-\sqrt{\lambda_l}\, t}a_{il}E_l(t).
$$

Hence 
$$
\langle Y_i(t),Y_j(t)\rangle=\langle \sum_{l=1}^{n-1}e^{-\sqrt{\lambda_l}\, t}a_{il}E_l(t),
\sum_{m=1}^{n-1}e^{-\sqrt{\lambda_m}\, t}a_{jm}E_m(t)\rangle
=\sum_{l=1}^{n-1}e^{-2\sqrt{\lambda_l}\, t}a_{il}a_{jl}.
$$
In order to estimate the determinant of the matrix
$$
A_k(t):=(\langle Y_i(t),Y_j(t)\rangle)_{1\leq i,j\leq k}
$$
 we consider the $k\times (n-1)$ matrices
$$
A:=(a_{il})\ \ \textup{and}\ \  B:=(a_{jl}e^{-2\sqrt{\lambda_l}\, t}).
$$
Then, by the Binet-Cauchy formula  (see
e.g. \cite{Gant}), we have 
$$
\det A_k(t)=\det (AB^{\top})=\sum_{1\leq s_1<s_2<\cdots<s_k\leq n-1}
\det A^{s_1\cdots s_k}\det B^{s_1\cdots s_k},
$$
where $A^{s_1\cdots s_k}$ denotes the $k\times k$ matrix consisting
of the columns of $A$ with indices $s_1\cdots s_k$ and similarly for $B$.
By definition $(a_{ij})_{1\leq i,j\leq n-1}$ is an orthogonal matrix and 
hence $|a_{ij}|\leq 1$.
By Proposition 1, the 
$\lambda_i$ are either zero or eigenvalues of the curvature operator given by $H_0$.
By assumption we have $k\geq\textup{rank}\ X$. Since the projection of the geodesics $\tau(x,t)$ to $S$ are regular by definition, at most $k-1$ of the Jacobi fields
$Y_i$ can be tangent to a flat in $X$. Thus,  for each matrix
$B^{s_1\cdots s_k}$, there is at least one column $s_i$, say, such that
 $\lambda_{s_i}=\alpha_{s_i}^2(H_0)>0$ for some positive root $\alpha_{s_i}$. Moreover, by the choice of $H_0$ as
 a barycentric direction, these values are \emph{uniformly bounded away
from zero} by $\lambda_*:=\min_{\alpha\in \Sigma^+}\alpha(H_0)>0$. We thus obtain the estimate
\begin{equation}\label{detA}
\det\, A_k(t)\leq C_1^2e^{-2\sqrt{\lambda_*}\, t},
\end{equation}
where $C_1^2$ depends only on $n$.
Inserting   (\ref{detA}) in (\ref{phi}) completes the proof of Lemma 1.
\hfill$\Box$

From Lemma 1 we immediately get (\ref{localest1}):
$$
\textup{vol}_{k+1}(\tilde{\Omega}_m^{r})=\int_{\Sigma_m^r}(\int_0^{\infty}\phi(x,t)dt)
\bar{\omega}(x)=C_2\textup{vol}_{k}({\Sigma}_m^{r}).
$$
\vspace{1ex}

{\sc STEP 4 (Volume of the cone over the singular part $\Sigma\setminus\Sigma^{r}$}):\ 

We claim that for $a\in [0,\infty)$ holds

\begin{equation*}
\textup{vol}_{k+1}(\tau(\Sigma\setminus\Sigma^{r}\times [0,a]))\leq 3^{k+1}
\textup{vol}_{k}({\Sigma}\setminus\Sigma^{r})\cdot a.
\end{equation*}
In fact, this immediately follows from 

\begin{lemma}  \label{1-lip}
The map $\tau:\Sigma\times  [0,\infty)\to X$ is 3-Lipschitz.
 \end{lemma}

\emph{Proof}.\ Let  $x,y\in \Sigma$. By construction the geodesics $\tau(x,t)$ and $\tau(y,t)$ are asymptotic. Thus by the convexity of the distance function in the CAT(0) space $X$
we have for $s,t\geq 0$ 
\begin{eqnarray*} d_{X}(\tau(x,s), \tau(y,t))\leq  d_{X}(\tau(x,s), \tau(y,s))+ d_{X}(\tau(y,s), \tau(y,t))\leq \\
\leq  d_{X}(\tau(x,0), \tau(y,0))+ \mid s-t\mid =d_{X}(x,y)+ \mid s-t\mid
\leq 3\sqrt{d_{X}(x,y)^2+ \mid s-t\mid^2}.
\end{eqnarray*}
\hfill$\Box$

\vspace{1ex}

{\sc STEP 5 (A chain $\Omega$ which fills  $\Sigma$)}:\ 
By Lemma 1 and since $I$ is finite we can find $t_0\geq 1$  so large that $\phi(x,t_0)\leq (\mid\!I\!\mid \,\textup{vol}_{k}({\Sigma}_m^{r}))^{-1}$ for all $m\in I$ (note that $I$ depends on $\Sigma$).

Consider a horosphere $\mathcal H$ based at $\omega\in X(\infty)$ (as in STEP 1) 
such that $\Sigma$ lies in the complement of the horoball $B_{\mathcal H}$ with boundary $\mathcal H$. Let $\pi_{\mathcal H}:X\to \mathcal H$
denote the horospherical projection. The image $\pi_{\mathcal H}(\Sigma)$ can  be written in the form
$$
\pi_{\mathcal H}(\Sigma)=\{\tau(x,t(x))\mid x\in \Sigma\}.
$$
We choose $\mathcal H$ in such a way that $t(x)\geq t_0$ for all $x\in \Sigma$.

By formula (\ref{dens}) and the choice of $t_0$  we the get
\begin{eqnarray*}
\textup{vol}_{k}(\pi_{\mathcal H}(\Sigma^{r}))=\sum_{m\in I}
\textup{vol}_{k}(\pi_{\mathcal H}(\Sigma_m^{r}))=
\sum_{m\in I}\int_{\Sigma_m^{r}}\phi(x,t(x))\bar{\omega}(x)\\ \leq
\sum_{m\in I}\int_{\Sigma_m^{r}}\phi(x,t_0)\bar{\omega}(x)
\leq \sum_{m\in I}(\mid\! I\!\mid\,\textup{vol}_{k}({\Sigma}_m^{r}))^{-1}\textup{vol}_{k}({\Sigma}_m^{r})=1.
\end{eqnarray*}

For the horosphere $\mathcal H$ chosen above we set $t_1:=\max\{t(x)\mid x\in \Sigma\}$ (we can assume that $t_1\geq 1$).
We now also choose the parameter $\epsilon$ of STEP 2 explicitely as $\epsilon:=(t_1\mid\!  I\! \mid)^{-1}$.

Notice that  the horospherical  projection $\pi_{\mathcal H}$  is
$1$-Lipschitz. By the definition of $\epsilon$ in STEP 2 and the explicit choice above we thus have $$\textup{vol}_{k}(\pi_{\mathcal H}({\Sigma}\setminus\Sigma^{r}))\leq \textup{vol}_{k}({\Sigma}\setminus\Sigma^{r})<\epsilon\mid\!I\!\mid=\frac{1}{t_1}\leq 1.$$
We conclude that
$$\textup{vol}_{k}(\pi_{\mathcal H}(\Sigma))\leq 1+1=2.$$ As this is 
 independent of $\Sigma$ and since $X$ has bounded geometry, there is a filling $\Omega_0$ of $\pi_{\mathcal H}(\Sigma)$ such that
 \begin{equation}\label{decomp}\textup{vol}_{k+1}(\Omega_0)\leq C_3,\end{equation}
 where  $C_3$ is again a uniform 
constant  depending only on $X$.
Finally, if $B_{\mathcal H}$ denotes the horoball in $X$ with boundary $\mathcal H$,  we set 
$$\Omega:=(\tilde{\Omega}\setminus B_{\mathcal H})\cup \Omega_0=(\tau(\Sigma\times \mathbb R_{\geq 0})\setminus B_{\mathcal H})\cup \Omega_0.$$
By construction  this $\Omega$  is a filling $(k+1)$-chain for
$\Sigma$.

\vspace{1ex}

{\sc STEP 6 (Volume  of $\Omega$)}:\ Let $t_0, t_1$ and $\epsilon=(t_1\mid\!I\!\mid)^{-1}$ be as in STEP 5.
From STEP 4 we have
\begin{eqnarray*}
\textup{vol}_{k+1}(\tau(\Sigma\setminus\Sigma^{r}\times [0,\infty))\setminus B_{\mathcal H})\leq \textup{vol}_{k+1}(\tau(\Sigma\setminus\Sigma^{r}\times [0,t_1])) \leq 
\\  \leq 3^{k+1}\textup{vol}_{k}({\Sigma}\setminus\Sigma^{r})\cdot t_1\leq 3^{k+1}\epsilon\mid\!I\!\mid t_1=3^{k+1}=:C_4
\end{eqnarray*}
This together with (\ref{localest}) and (\ref{decomp}) yields
\begin{eqnarray*}
\textup{vol}_{k+1}({\Omega})=
\textup{vol}_{k+1}(\tau(\Sigma^{r}\times [0,\infty))\setminus B_{\mathcal H})+
\textup{vol}_{k+1}(\tau(\Sigma\setminus\Sigma^{r}\times [0,\infty))\setminus B_{\mathcal H})+
\textup{vol}_{k+1}(\Omega_0)\\ \leq\textup{vol}_{k+1}(\tilde{\Omega}^{r})+C_4+C_3\leq 
C_2\textup{vol}_{k}(\Sigma^{r})+C_4+C_3\leq C_2\textup{vol}_{k}(\Sigma)+C_3+C_4\prec \textup{vol}_{k}(\Sigma).
\end{eqnarray*}

This estimate completes the proof of of the linear upper bound in Theorem 1 (ii).

\subsubsection{Linear lower bounds}

In order to establish  linear lower bounds for the filling functions in Theorem 1 (ii), it suffices to provide
arbitrarily large $k$-cycles $\Sigma$, with the property that any $(k+1)$-filling $\Omega$ satifies
$\textup{vol}_{k+1}( \Omega) \succ \textup{vol}_{k}(\Sigma)$. We construct such cycles 
via the following  lemma.
Recall that every symmetric space $S$ of noncompact type contains (many) copies of 
real hyperbolic planes
$\mathbb H ^2$ (see \cite{Hel}, IX.2). 

\begin{lemma} Let $S$ be an $n$-dimensional symmetric space of noncompact type  
 with base point $x_0$. Let $\mathbb H ^2\subset S$ be a totally geodesic 
hyperbolic plane containing $x_0$. For any $2\leq k\leq n$ there exists a
  $(k-2)$-dimensional subspace $\mathfrak{q}$ of $\mathfrak{p}\cong T_{x_0}S$ orthogonal
 to $\mathfrak{h}\cong T_{x_0}\mathbb H ^2$ such that
the $k$-dimensional submanifold
$$
W:= \exp (\mathfrak{h}\oplus\mathfrak{q})\cdot x_0,
$$
of $S$, which is diffeomeorphic to 
$\mathbb H ^2\times \exp\mathfrak{q}\cdot x_0=:\mathbb H ^2\times Q$, has the following property:

Let $\Sigma$ be the topological sphere obtained as the intersection of $W$ with the sphere $S_R(x_0)$ of radius $R>>1$
 and center $x_0$ in $S$. 
Then the $k$-dimensional Hausdorff measure (or Riemanian volume) of any (smooth)
 filling $\Omega$  of $\Sigma$ \underline{in $S$} satisfies 
$$\textup{vol}_{k}( \Omega) \succ e^R\succ \textup{vol}_{k-1}(\Sigma).$$
 \end{lemma}

\emph{Proof}.\  The proof  is by induction on the dimension of $W$ and
 uses the  coarea formula.

We start the induction with $k=2$. 
In this case we take $\mathfrak{q}:=\{0\}$, i.e., $W=\mathbb H^2$.
Consider the circle $\Sigma_1:=S_R(x_0)\cap W\subset W$. Fill $\Sigma_1$ with
 any $2$-chain $\Omega_2$ \underline{in $S$}.
The orthogonal projection  of $S$ to the totally geodesic submanifold $W=\mathbb H^2$ is $1$-Lipschitz. Hence, if $\Omega_2'$ is the orthogonal projection of $\Omega_2$ to $W$, we have $\textup{vol}_2(\Omega_2)\geq \textup{vol}_2(\Omega_2')$, i.e, the filling $\Omega_2'$ of $\Sigma_1$ \underline{in $W=\mathbb H^2$} is even smaller than the filling in $S$. The minimal filling of the circle $\Sigma_1$ in $W=\mathbb H^2$ is the disc $D^2_R$ of radius $R$. Hence
$$ \textup{vol}_{2}( \Omega_2')\geq\textup{vol}_{2}(D^2_R)\sim e^R\sim \textup{vol}_{1}(\Sigma_1)$$
and the claim follows in this case.

In order to continue the proof by induction, we
assume that the lemma is proved for suitable submanifolds $W^2, W^3,\ldots W^k$ with $\dim W^i=i$.
We now  wish to prove the lemma for a $W^{k+1}$ of dimension $k+1$. To that end we choose a geodesic $c$ in $S$ orthogonal to $W^k$ and define $W^{k+1}$ as the union of the parallel translates of $W^k$ along $c$. Let  $\pi:S\to c$ denote  the orthogonal projection of $S$ to the geodesic $c$.  Since $\pi$ is $1$-Lipschitz,  Eilenberg's coarea formula (see \cite{BZ}, 13.3) yields
for any filling $\Omega_{k+1}$ of $\Sigma_k:=S_R(x_0)\cap W^{k+1}$ in $S$
\begin{equation}\label{ind-co}
\textup{vol}_{k+1}(\Omega_{k+1})\succ  \int^{\infty}_{-\infty}\textup{vol}_k(\pi^{-1}(t)\cap \Omega_{k+1})\, dt\geq \int^{+1}_{-1}\textup{vol}_k(\pi^{-1}(t)\cap \Omega_{k+1})\, dt.
\end{equation}

The last integral is geater than the minimal filling volume of the intersection of $S_R(x_0)$ with the ``thickening'' of $W^k=\pi^{-1}(0)\cap W^{k+1}$ given
by
$$
\bigcup_{t\in (-1,+1)}\pi^{-1}(t)\cap W^{k+1}.
$$
The latter consists of copies of $W^k$ (note that
$\pi$ commutes with parallel translation along $c$).

By Sard's theorem almost every $t\in[-1,1]$ is a regular value of $\pi\mid \Omega_{k+1}$ and 
$\pi\mid \Sigma_{k}$. Hence 
$$
\partial(\pi^{-1}(t)\cap\Omega_{k+1})=\pi^{-1}(t)\cap \partial\Omega_{k+1}=\pi^{-1}(t)\cap \Sigma_{k},
$$
i.e., 
$\pi^{-1}(t)\cap\Omega_{k+1}$ is a filling of $\pi^{-1}(t)\cap \Sigma_{k}$ for a.e. $t\in [-1,1]$  (see e.g. \cite{GP} 2.1.). 
Now we have
$$
\pi^{-1}(0)\cap \Sigma_k=\pi^{-1}(0)\cap W^{k+1}\cap S_R(x_0)= W^k\cap S_R(x_0)=\Sigma_{k-1}.$$
We thus can apply the induction hypothesis and get that $\pi^{-1}(0)\cap \Omega_{k+1}$ and hence also $\pi^{-1}(t)\cap \Omega_{k+1}$ for $t\in [-\frac{1}{100},\frac{1}{100}]$, say,  has Hausdorff measure (or Riemannian volume) $$\textup{vol}_{k}(\pi^{-1}(t)\cap\Omega_{k+1})\succ e^R.$$
Inequality (\ref{ind-co}) then yields
 $\textup{vol}_{k+1}(\Omega_{k+1})\succ e^R$.

Finally, we need to show that $\textup{vol}_k(\Sigma_k)\prec e^R$. But $\Sigma_k$ is a sphere of radius $R$ in $W^{k+1}$ with center $x_0$ and  $W^{k+1}$ contains a hyperbolic plane through $x_0$. The claim then follows e.g. from the computation of volumes of spheres via Jacobi fields (see \cite{Sak}, II.5, Lemma 5.4). This finishes the induction proof of Lemma 2.
\hfill$\Box$

\vspace{1ex}

We proceed to prove  the lower bounds in Theorem 1 (ii).
Let $\Sigma_k$ be a $k$-sphere of radius $R$ in  $S\subset X=S\times E\times B$ as in Lemma 2 and let $\Omega_{k+1}$ be an arbitrary filling of $\Sigma_k$ in $X$. Let $\pi:X\to S$ be the orthogonal projection. As $S$ is a convex subspace
of $X$ this again is a 1-Lipschitz map. Hence, by Lemma 2,
$$
\textup{vol}_{k+1}(\Omega_{k+1})\geq  \textup{vol}_{k+1}(\pi(\Omega_{k+1}))\succ  \textup{vol}_{k}(\Sigma_k).
$$

\vspace{1ex}

In order to get 
$k$-cycles having at least linear fillings for any chain in $X$ for $k$ up to $\dim X-1=\dim S+\dim E+\dim B-1$ we choose a flat $F$ in $E\times B$ 
and consider the convex subspace $S\times F$ in $X$. Note that this actually is a manifold of dimension equal to $\dim X$.
For the construction of the  $\Sigma_k$  with $k\geq\dim S$ we then proceed exactly as in the last step
of the proof of the Lemma 2 by inductively defing suitable $W^{k+1}\subset S\times F\subset X$ using  successively parallel translation along and orthogonal projection to pairwise orthogonal geodesics in $F$. Finally, since the orthogonal projection from $X$ to $S\times F$ is $1$-Lipschitz, we get the linear filling inequality by the same argument as above.

\section{The proof of Theorem 2}

We have the following general result relating filling and Dehn functions 

\begin{proposition}  Let $X_1$ and $X_2$ be two $k$-connected manifolds or simplicial complexes on which a group $\Gamma$ acts geometrically, then
$$
\delta^k_{X_1}\sim \delta^k_{X_2}\  \ \ \textup{and}\ \ \ FV_{X_1}^{k+1}\sim FV_{X_2}^{k+1}.
$$
Moreover $\delta^2_{X_1} \prec FV_{X_1}^{3}$  and $\delta^k_{X_1} \sim FV_{X_1}^{k+1}$ for $k\geq 3$.
\end{proposition}

\emph{Proof}.\ The first  assertion is proved in \cite{AWP}, 6. Corollary 3, for a (simplicial) version of $\delta^k$. This is equivalent to the above (compare \cite{Yo1}).
The second claim follows from \cite{Gr1}, 2.A', for upper bounds and  \cite{BBFS}, 2.6.(4),
for lower bounds.
\hfill$\Box$

\vspace{2ex}

We emphasize that Young constructed  examples of  groups $\Gamma$ acting geometrically on a space $X$ for which 
$\delta^2_{\Gamma} \not\sim FV_{X}^{3}$ (see \cite{Yo1}, Corollary 6). If, however,  $X=S\times E\times B$ is as in Theorem 1, then we also have $ \delta^2_{X} \succ  FV_{X}^{3}$. In fact,
by Theorem 1 and the explicit construction of 
lower bounds by filling spheres by balls in sections 2.1 and 2.2.2 we have
\begin{eqnarray*}
 \delta^2_{X}(l) \succ l^{\frac{3}{2}} \sim FV_{X}^{3}(l) & \textup{ if } \textup{rank}\ X\geq 3\\
\delta^2_{X}(l) \succ l \sim FV_{X}^{3}(l) & \textup{ if } \textup{rank}\ X=1 \ \textup{or } 2.
\end{eqnarray*}
Hence in that case there are no cycles that are ``harder to fill'' than spheres. Together with this observation and the fact (by \cite{AWP}) that $\delta^k_{\Gamma}$ is a quasi-isometry invariant, i.e., $\delta^k_{\Gamma}\sim \delta^k_{X}$,
 Theorem 2 follows immediately from Theorem 1 and Proposition 3.

\section{Basic Examples}

Examples of groups $\Gamma$ to which Theorem 2 applies include:
\begin{enumerate}
 \item
 Crystallographic groups, i.e.,
discrete groups of isometries  acting cocompactly on some Euclidean space $\mathbb E^n$
(see e.g. \cite{VS}, Part II, Ch. 3).
\item
Uniform lattices in  semisimple real Lie groups $G$, which act cocompactly on the assocated symmetric 
space $S=G/K$ of noncompact type (see e.g. \cite{Bor}).
\item Uniform $S$-arithmetic groups, which act
cocompactly on products of  symmetric spaces of noncompact type  and 
 of Bruhat-Tits buildings. We briefly describe these last examples and some of their properties in the next paragraph.
\end{enumerate}

\subsection{S-arithmetic groups}

We briefly recall the definition of $S$-arithmetic groups.
For more details and informations we refer to \cite{BS} and  \cite{Mar}, Chapter I.3.   

 A global  field $\mathbb K$ is either an algebraic number field, i.e., a finite extension
of $\mathbb Q$, or the function field of an algebraic curve over a finite field, i.e., a finite extension
of $\mathbb F_q(T)$, the field of rational functions in one variable over the finite field with $q$ elements.
In the first case the characteristic of $\mathbb K$ is zero and in the second it is positive. Completions of global fields with respect to archimedean or non-archimedean valuations are local fields.

 Let $\mathbb K$ be a global field and let $S$ be a finite set of inequivalent valuations of $\mathbb K$ including all archimedean ones. The ring of $S$-integers of $\mathbb K$ is defined as $$\mathcal O_S:=\{k\in \mathbb K\mid \ \mid\!k\!\mid_v\leq 1\ \textup{ for all } v\notin S\}.$$ Let  $\mathbf G\subset \mathbf{GL}_N$ be a
connected, reductive  $\mathbb K$-group. The associated S-arithmetic group is the group of $S$-integers, $\mathbf G(\mathcal O_S)=\mathbf G\cap\mathbf{GL}_N(\mathcal O_S)$. Its diagonal embedding  into the locally
compact group $\mathbf G_S:=\prod_{v\in S}\mathbf G({\mathbb K}_v)$  (here ${\mathbb K}_v$ is the completion of $\mathbb K$ with respect to $v$) is an irreducible lattice. The reductive group $\mathbf G_S$ in turn acts isometrically
on an associated  CAT(0) space $X$ as in Theorem 1 or 2, i.e., the metric product  of a symmetric space of noncompact type, a Euclidean space  and a product of Euclidean (Bruhat-Tits) buildings. 

\subsection{Some (non-)existence results}

We list a few  results concering existence 
and non-existence of uniform $S$-arithmetic groups.

\begin{itemize}
\item (Borel-Harder, \cite{BH})\  Let $\mathbf G$ be a reductive $\mathbb L$-group over a local field $\mathbb L$ of characteristic zero. Then the group of $\mathbb L$-points, $\mathbf G(\mathbb L)$, contains uniform lattices. 
 \item (Borel-Harish-Chandra, Margulis, \cite{Mar}, I.3.2.7)\ Let $\mathbb K$ be a global field.
If $\mathbf G$ is  \emph{$\mathbb K$-anisotropic} and if one  chooses  the finite set $S$ in such a way
 that $\mathbf G$ is \emph{${\mathbb K}_v$-isotropic} for any $v\in S$, then (the diagonal embedding of) $\mathbf G(\mathcal O_S)$ is a uniform lattice in $\mathbf G_S$  and thus acts cocompactly on $X$. Note that many such $S$ exist, since the set of valuations $v$ such that
$\mathbf G$ is anisotropic over ${\mathbb K}_v$ is finite (see \cite{Spr}, Lemma 4.9).
\item (Harder, \cite{Ha}, 3. Kor. 1)\ Only groups of type $\mathbb A$ have anisotropic forms over global fields of positive characteristic. In particular, if at least one of the factors 
$\mathbf G({\mathbb K}_v)$  in $\mathbf G_S$ has $\textup{char}\ {\mathbb K}_v>0$ and is not of type $\mathbb A$, then $\mathbf G_S$ cannot contain uniform lattices.
\item If $\mathbf K$ is a global field and $\mathbf G$ is $\mathbb K$-isotropic, then $\mathbf G(\mathcal O_S)$ is a non-uniform lattice. 
\end{itemize}

\subsection{Finiteness properties}

The very definition of higher-dimensional Dehn functions for some group $\Gamma$ requires that $\Gamma$  acts on highly connected spaces: $k$-dimensional Dehn functions are meaningful only for groups of finiteness type $F_{k+1}$.

In the case of \emph{number fields} we have the following result:
\begin{quote}
(Raghunathan \cite{Ra}, Borel-Serre \cite{BS})\ \ Any $S$-arithmetic subgroup $\mathbf G(\mathcal O_S)$ (uniform or not) of a reductive group $\mathbf G$
defined over a global number field is of finiteness type $F_{\infty}$ 
\end{quote}

In the case of \emph{function fields} there are restrictions:

\begin{quote} (Serre \cite{Se}, Kropholler-Mislin \cite{KM}, Bux-Wortman \cite{BW})\ \
An $S$-arithmetic subgroup $\mathbf G(\mathcal O_S)$ of a reductive group $\mathbf G$ defined over a function field $\mathbb K$ is of type $F_{\infty}$ if and only if $\mathbf G$ is $\mathbb K$-anisotropic (or, equivalently, if $\mathbf G(\mathcal O_S)$ is uniform).
\end{quote}

Hence  for \emph{ uniform} $S$-arithmetic groups Dehn functions are always defined in any dimension. 
By the above result of Borel and Serre this also holds for non-uniform $S$-arithmetic groups in the number field case.  We emphasize however that this is no longer true for non-uniform $S$-arithmetic subgroups of  $\mathbb K$-isotropic groups $\mathbf G$ defined over \emph{function fields}.
Recall that the finiteness length of a group $\Gamma$ is the maximal $m$ such that $\Gamma$ has finiteness type $F_{m}$.  In \cite{BW} and \cite{BGW} it is proved
  that if $\mathbf G(\mathcal O_S)$
is non-uniform (and $\mathbf G$ defined over a function field), then its finiteness length is $r-1$, where $r$ is the sum of the local ranks (or, equivalently, the Euclidean rank of the associated  product of Euclidean buildings on which $\mathbf G(\mathcal O_S)$ acts).
In particular, if $k\geq r-1$, the
the $k$-dimensional  Dehn function is  not defined for such groups. We adress a conjectural
picture about Dehn functions for non-uniform lattices in the next section.

\section{Some conjectures and results for non-uniform lattices}

\subsection{A general conjecture}
If a group $\Gamma$ acts geometrically on a metric space $X$, then $\Gamma$ is quasi-isometric to $X$.
Examples are {\it uniform} S-arithmetic groups acting on $X=S\times E\times  B$ (see Section 4).
 This is no longer true for 
 {\it non-uniform} S-arithmetic groups like $SL_n(\mathbb Z[\frac{1}{p}])$ which in fact can be strongly distorted and thus make the investigation of filling problems
considerably more subtle. We  posit the  following conjectural picture
 (compare also the discussion in Gromov's book \cite{Gr2}, 5D(c)).

\begin{conjecture} Let $\Gamma$ be an $S$-arithmetic subgroup of a reductive algebraic group defined  over a global  field $\mathbb K$ such that $\Gamma$ acts as an (irreducible) \emph{non-uniform} lattice on $X=S\times E\times  B$, the associated metric product  of symmetric spaces,
Euclidean spaces  and Euclidean buildings. 
Let $r$ be the Euclidean rank of $X$ and $n$ the dimension of $X$. 
 Then the higher-dimensional Dehn functions of $\Gamma$ satisfy:$$
\textup{(i)}\ \hspace{1cm} \delta_{\Gamma}^k(l)\sim l^{\frac{k+1}{k}},\ \ \ \ \ \textup{ if}\ \  \ k\leq r-2.
$$
Moreover, if $\mathbb K$ has characteristic zero (i.e., if $\mathbb K$ is a number field), then

\begin{eqnarray*}
\textup{(ii)} & \  \delta_{\Gamma}^{r-1}(l)\sim \exp l; &\\
\textup{(iii)} & \  \delta_{\Gamma}^k(l)\sim l^{q(k)},& \!\!\!\textup{with }   1\leq q(k)\in\mathbb Q, \  \
\textup{ if}\ \ r\leq k\leq n-2;\\
\textup{(iv)} &   \delta_{\Gamma}^{n-1}(l)\sim  l.
\end{eqnarray*}
\end{conjecture}

\vspace{1ex}

\subsection{Some evidence for Conjecture 1}

In the following we discuss partial results which provide some evidence for Conjecture 1.

We will use the elementary
\begin{lemma} \label{convlip} Let $X$ be a CAT(0) space and $X_0$ a length subspace. Suppose that 
$C$ is a complete convex subspace both of $X$ and $X_0$ in the induced metrics. Then the restriction of the orthogonal projection
$\pi:X\to C$ to $X_0$ is $1$-Lipschitz.
 \end{lemma}
\emph{Proof}.\ By \cite{BH} Ch.II, Prop. 2.4(4), the projection $\pi$  is $1$-Lipschitz. Moreover, since $C$ is convex in $X_0$ as well as in $X$, we have for all $x,y\in X_0$
$$
\ \ d_{X_0}\!\mid_C(\pi(x),\pi(y))=d_C(\pi(x),\pi(y))=d_{X}\!\mid_C(\pi(x),\pi(y))\leq d_X(x,y)\leq d_{X_0}(x,y). \hspace{1.7cm}\Box
$$

\subsubsection{Linear Dehn functions in top dimension}

The next theorem confirms  Conjecture 1(iv) for lattices in semisimple groups.

\begin{theorem} \label{topdim} Let  $X=G/K$ be an $n$-dimensional symmetric space of noncompact type
(without Euclidean factor) and let $\Gamma\subset G$ be a lattice. Then  the top-dimensional Dehn function of $\Gamma$ is linear:
$$\delta_{\Gamma}^{n-1}(l)\sim l.$$
 \end{theorem}

\emph{Proof}.\ If $\Gamma$ is uniform the claim follows from Theorem 2. If $\Gamma$ is non-uniform, it does not act cocompactly on $X$, but there is a suitable subspace $X_0$ of $X$ obtained by deleting a $\Gamma$-invariant family of horoballs on which $\Gamma$ acts cocompactly
(see \cite{GR} for $\textup{rank}\, X= 1$ and \cite{Le} for $\textup{rank}\, X\geq 2$). Since $\Gamma$ is not virtually solvable (see \cite{OV}, Ch. 4, Thm. 3.6)  it contains a free group by the Tits alternative and thus is non-amenable. Moreover, $X$ and  $X_0$ have bounded geometry. Thus, by \cite{Gr2} $0.5.A_6$, the $n$-dimensional (Hausdorff) volume of any domain 
$\Omega\subset X_0$ with $(n-1)$-dimensional measurable boundary $\partial\Omega$ satisfies
a linear isoperimetric inequality
$$
\textup{vol}_n(\Omega)\leq \textup{const. }\textup{vol}_{n-1}(\partial\Omega)
$$
for  a uniform constant depending only on $X_0$.
This observation yields a linear upper bound for the $n$-dimensional filling function (and hence 
the $(n-1)$-dimensional Dehn function).

In order to get a linear lower bound for the $n$-dimensional filling function, we choose a geodesic $c$ in $X_0$ that is mapped to a closed geodesic in $X/\Gamma$. Such a $c$ exist: One can e.g. choose $c$ in a $\Gamma$-closed flat (see \cite{Mo}, Lemma 8.3').
Then take (a small) $\delta>0$ such that the $\delta$-tube around 
$c$ is completely contained in $X_0$ (replace $X_0$ by a $2\delta$-neighbourhood of $X_0$ if necessary). The $\delta$-tube $U_{\delta}(c_l)$ of a segment of $c$ of length $l$ is a convex subset of $X$ and also of $X_0$ (see \cite{BH} Ch.II, Cor. 2.5). By Lemma \ref{convlip} the  orthogonal projection from $X$ to $U_{\delta}(c_l)$ restricted to $X_0$ is 1-Lipschitz. 
Hence any filling of $\partial U_{\delta}(c_l)$ in $X_0$ has at least the same volume as a filling 
in $U_{\delta}(c_l)$. Since $X$ and $X_0$ have bounded geometry, we have  $$\textup{vol}_n(U_{\delta}(c_l))\sim l\sim \textup{vol}_{n-1}(\partial U_{\delta}(c_l))
$$ 
which yields the claimed linear lower bound for the $n$-dimensional filling  function.
\hfill$\Box$

\subsubsection{Proof of Conjecture 1 for real and complex hyperbolic lattices}
Here we  confirm Conjecture 1 for non-uniform lattices (arithmetic or not) in $SO(n,1)$ resp. $SU(n,1)$, the group of orientation preserving isometries  of real hyperbolic space $H^n\mathbb R$ resp. complex hyperbolic space $H^n\mathbb C$. Notice that in this case the Euclidean rank of 
$X=H^n\mathbb R$ (resp. $H^n\mathbb C$) is one, so that Conjecture 1 reduces to parts (iii) and (iv).
Note further that  \emph{uniform} lattices in $SO(n,1)$ resp. $SU(n,1)$ are covered by Theorem 2: As the rank is one, all higher-dimensional Dehn functions are linear.

\begin{theorem}\label{hyperbol} Let $\Gamma$ be a non-uniform lattice in $SO(n,1)$.
Then the  higher-dimensional Dehn functions of $\Gamma$ satisfy:
\begin{eqnarray*}
\textup{(a)} &\ \ \delta_{\Gamma}^k(l)\sim l^{\frac{k+1}{k}}&\ \ \ \textup{for}\ \ 1\leq k\leq n-2,\\
\textup{(b)} & \ \ \delta_{\Gamma}^{n-1}(l)\sim l.&\\
\end{eqnarray*}
 \end{theorem}

\emph{Proof}.\ It is well known that $\Gamma$ acts geometrically on a CAT(0) space $X_0$ obtained from
real hyperbolic space $H^n\mathbb R$ by removing a $\Gamma$-invariant family of disjoint horoballs
(see  \cite{BH}, Thm. 11.27, Cor. 11.28). The Euclidean rank of $X_0$ is $n-1$ and (isolated) maximal flats
are provided by the boundary horospheres of the deleted horoballs. By Wenger's result (Proposition \ref{cat0-iso}) the simply connected CAT(0) space $X_0$ satisfies 
Euclidean isoperimetric inequalities below the rank:
$$ \textup{FV}_{X_0}^{k+1}(l)\prec l^{\frac{k+1}{k}}\ \ \ \textup{for}\ \ 1\leq k\leq n-2,$$
Lower Euclidean bounds are obtained as follows. Take a maximal flat in $X_0$, i.e., a boundary horosphere $H$,
and a round $k$-sphere $\Sigma_k\subset H$. By Lemma \ref{convlip} the orthogonal projection from $X_0$ to the convex subspace $H$ is $1$-Lipschitz. Hence the projection of any minimal filling $\Omega_{k+1}$
of $\Sigma_k$ in $X_0$ orthogonally projects to an even smaller minimal  filling $\Omega_{k+1}'$ of $\Sigma_k$ in $H$ and hence is a $(k+1)$-ball, which yields the Euclidean lower bound:
$$
\textup{vol}_{k+1}(\Omega_{k+1})\geq \textup{vol}_{k+1}(\Omega_{k+1}')\sim (\textup{vol}_{k}(\Sigma_{k}))^{\frac{k+1}{k}}.
$$

Claim (b) on the top-dimensional Dehn function follows from Theorem \ref{topdim}.
\hfill$\Box$
\begin{theorem}\label{hypcomp} Let $\Gamma$ be a non-uniform lattice in $SU(n,1)$.
Then the  higher-dimensional Dehn functions of $\Gamma$ satisfy:
\begin{eqnarray*}
\textup{(a)} &\ \ \delta_{\Gamma}^k(l)\sim l^{\frac{k+1}{k}}&\ \ \ \textup{for}\ \ 1\leq k\leq n-2,\\
\textup{(b)} &\ \ \delta_{\Gamma}^{n-1}(l)\sim l^{\frac{n+1}{n-1}},&\\
\textup{(c)} &\ \ \delta_{\Gamma}^k(l)\sim l^{\frac{k+2}{k+1}}&\ \ \ \textup{for}\ \ n\leq k\leq 2n-2,\\
\textup{(d)} & \ \ \delta_{\Gamma}^{2n-1}(l)\sim l.& \\
\end{eqnarray*}
 \end{theorem}

\emph{Proof}.\ As in the proof of Theorem \ref{hyperbol} we use the fact that $\Gamma$ acts geometrically on a  space $X_0$ obtained from
real hyperbolic space $H^n\mathbb C$ by removing a $\Gamma$-invariant family of disjoint horoballs
(see e.g. \cite{GR} or \cite{Le}). Notice that in contrast to the real hyperbolic case, this $X_0$ is {\it not} a CAT(0) space.
The boundary horospheres of $X_0$ can be identified with $(2n-1)$-dimensional Heisenberg groups $H_{2n-1}$.
Moreover, given a boundary horosphere ${\mathcal H}\subset \partial X_0$, there is a $1$-Lipschitz
retraction $r:X_0\to \mathcal H$ (see \cite{Le}). Hence, if $\Sigma$ is a $k$-cycle in $\mathcal H$ any filling $\Omega$ of $\Sigma$
\underline{in $X_0$} cannot be of smaller volume than that of the minimal filling of $\Sigma$ \underline{in $\mathcal H$}. This yields lower bounds for the filling functions of $X_0$:
$$
\textup{FV}_{X_0}^{k+1}\succ \textup{FV}_{\mathcal H}^{k+1}\ \ \textup{ for }\ 1\leq k\leq 2n-2.
$$
In order to obtain upper bounds we argue as in \cite{Gr2} 5.D.5(c). Given a $k$-cycle $\Sigma$ in $X_0$ we first take a filling $\widetilde{\Omega}$ of $\Sigma$ \underline{in $X$} with $\textup{vol}_{k+1}(\widetilde{\Omega})\prec
\textup{vol}_{k}(\Sigma)$ (such a linear filling exists by Theorem 1). Then we consider the intersection $\widetilde{\Sigma}:=\widetilde{\Omega}\cap \partial X_0$. This is a $k$-cycle in $\partial X_0$. 
We fill $\widetilde{\Sigma}$  by a minimal $(k+1)$-chain $\widehat{\Omega}$ \underline{in  $\partial X_0$}.
Together with $\widetilde{\Omega}\cap X_0$ this eventually yields a $(k+1)$-chain $\Omega$ which fills $\Sigma$
\underline{in $X_0$}.

Now assume that $\textup{vol}_{k}(\Sigma)=l$. By taking a $\epsilon$-neighbourhood of $\widetilde{\Sigma}$ we see that
$$
\textup{vol}_{k}(\widetilde{\Sigma})\sim 
\textup{vol}_{k+1}(U_{\epsilon}(\widetilde{\Sigma})\cap \widetilde{\Omega})<\textup{vol}_{k+1}(\widetilde{\Omega})
\prec\textup{vol}_{k}(\Sigma).
$$ Hence also have $\textup{vol}_{k}(\widetilde{\Sigma})\prec l$. This yields
\begin{eqnarray*}
\textup{vol}_{k+1}(\Omega)\leq \textup{vol}_{k+1}(\widetilde{\Omega}\cap X_0)+
\textup{vol}_{k+1}(\widehat{\Omega})
\leq \textup{vol}_{k+1}(\widetilde{\Omega})+
\textup{vol}_{k+1}(\widehat{\Omega})\leq l + \textup{FV}^{k+1}_{\partial X_0}(l),
\end{eqnarray*}
where the last estimate follows from the fact that $\widehat{\Omega}$ is a minimal filling in $\partial X_0$ of the cycle
$\widetilde{\Sigma}$, which by the above also has $k$-volume $\leq l$. Since $\partial X_0$ is the union of disjoint horospheres we further have $\textup{FV}^{k+1}_{\partial X_0}\sim \textup{FV}^{k+1}_{\mathcal H}$.
We conclude that
$$\textup{FV}^{k+1}_{X_0}(l)\leq l + \textup{FV}^{k+1}_{\mathcal H}(l)\sim l +\textup{FV}^{k+1}_{H_{2n-1}}(l).
$$
The higher-dimensional filling functions of Heisenberg groups have been determined by Young in \cite{Yo2}.
As they are superlinear we obtain that
$$\textup{FV}^{k+1}_{X_0}\prec\textup{FV}^{k+1}_{H_{2n-1}}\sim  \textup{FV}^{k+1}_{\mathcal H}.$$
Together with the lower bound above this eventually yields for the filling functions of $X_0$:
$$\textup{FV}^{k+1}_{X_0}\sim \textup{FV}^{k+1}_{\mathcal H}\sim \textup{FV}^{k+1}_{H_{2n-1}}\ \ \textup{ for }\ 1\leq k\leq 2n-2.$$
As $\Gamma$ is quasi-isometric to $X_0$  we get for the Dehn functions for $\Gamma$:
$$
\delta_{\Gamma}^k\sim \delta_{X_0}^{k}\sim \delta_{\mathcal H}^{k}\sim \delta_{H_{2n-1}}^{k}\ \ \textup{ for }\ 1\leq k\leq 2n-2.
$$

Parts (a), (b) and (c) of the theorem then follow  from  \cite{Yo2}.

Claim (d) on the linear top-dimensional Dehn function follows from Theorem \ref{topdim}.
\hfill$\Box$

\vspace{1ex}

\subsubsection{Low dimensional Euclidean lower bounds}
The last argument of the proof of Theorem \ref{hyperbol} (a) generalizes to non-uniform lattices of semisimple Lie groups:

\begin{theorem}\label{euclower} Let  $\Gamma$ be an irreducible, non-uniform lattice in a semisimple Lie group $G$
acting on a symmetric space $X=G/K$ of rank $r$. Then the Dehn functions of dimension between $1$ and $r-2$  have  Euclidean lower bounds:
  $$\delta_{\Gamma}^k(l)\succ l^{\frac{k+1}{k}} \ \ \textup{if}\ \ 1\leq  k\leq r-2.$$
 \end{theorem}

\emph{Proof}.\  A theorem of Mostow asserts that there are many closed maximal flats in $X\backslash \Gamma$, see
\cite{Mo}, Lemma 8.3, 8.3'.
Pick one such flat and its universal cover $F$ in $X_0$ (a length subspace of $X$ quasi-isometric to $\Gamma$, see \cite{GR} and 
\cite{Le}).   Then the restriction of the orthogonal projection  $\pi:X\to F$ to $X_0$ is $1$-Lipschitz
by Lemma \ref{convlip}
 and the claim follows from the solution of the Euclidean isoperimetric problem in $F$.\hfill$\Box$

\subsubsection{Exponential (r-1)-dimensional Dehn functions}

Conjecture 1(ii), i.e., that the $(r-1)$-dimensional Dehn function is exponential,  has been proved in a number of cases.

 By a result of Gromov (see \cite{Le} for a
detailed proof) Dehn functions (of any dimension) of linear groups are at most exponential.

Wortman, in  \cite{Wo1}, established  lower exponential bounds  for many non-uniform arithmetic subgroups of semisimple algebraic $\mathbb Q$-group.
The cases where the conjecture for such groups is still open are those of  relative $\mathbb Q$-type
$BC_n, G_2, F_4$ and $E_8$. 

In \cite{LP} we proved exponential lower bounds for the 1-dimensional Dehn function of (irreducible) non-uniform lattices in all semisimple Lie groups of rank $2$.

By Taback \cite{Ta}, the $1$-dimensional Dehn function for the (simplest) $S$-arithmetic group
$PSL_2(\mathbb Z[\frac{1}{p}])$, $p$ a prime number, is exponential.

\subsection{The Bux-Wortman conjecture}

A  unified approach to large scale properties of $S$-arithmetic groups on function fields or number fields  has been formulated in a conjecture by
 Bux and Wortman \cite{BW}. 
It is based on a geometric version of reduction theory and uses the concept of coarse manifolds to measure the distorsion of the metrics
in various dimensions. 

A \emph{coarse $n$-manifold of scale $s$} in a metric space $X$
is the image of a map $C:M\to X$ from a triangulated $n$-dimensional manifold   $M$ (with boundary) such that any two adjacent vertices of the triangulation are 
mapped to points of distance at most $s$ in $X$.
The boundary of a coarse manifold is the image under $C$ of the boundary of $M$.  The \emph{volume of a coarse manifold} is defined as the number
of vertices in the triangulation.  
A subset $Y\subset X$ is called \emph{(homotopically) undistorted in dimension $k$}, if there are
positive constants $s_1$ and $s_2$ and a linear function $L:\mathbb R\to\mathbb R$ such that the following holds:

\begin{quote}
For any coarse $k$-ball $B$ in $X$ of scale at most $s_1$ whose boundary is a coarse sphere in $Y$,
there is a coarse $k$-ball $B'$ in $Y$ of scale at most $s_2$ with the same boundary, such that $\textup{vol}_k(B')\leq L(\textup{vol}_k(B))$.
\end{quote}

The \emph{(homotopical) distorsion dimension} of $Y$ in $X$ is the largest dimension $k$ such that $Y$ is homotopically undistorted in dimension $k$.

Consider now a reductive $\mathbb K$-group $\mathbf G$ defined over a global field $\mathbb K$ and a finite set $S$ of valuations  as in Section 4.1.
Let $\mathbf G_S$ be the semisimple Lie group $\mathbf G_S=\prod_{v\in S}\mathbf G(\mathbb K_v)$
endowed with a left-invariant metric. Notice that this metric space is quasi-isometric to the product 
$X=S\times B$ of a symmetric space and Euclidean buildings on which $\mathbf G_S$ acts isometrically. The metric on $\mathbf G_S$ restricts to the diagonally embedded lattice
$\mathbf G(\mathcal O_S)$. The Euclidean rank $r$ of $\mathbf G_S$ (or $X$) is the sum of the local ranks $r:=\sum_{s\in S}\textup{rank}_{\mathbb K_v}\mathbf G$.

\begin{conjecture} \textup{(Bux-Wortman, \cite{BW})}\ The $S$-arithmetic subgroup $\mathbf G(\mathcal O_S)$ has distortion dimension $r-1$ as a subgroup of the reductive group $\mathbf G_S$ with respect to the metrics 
defined above. 
\end{conjecture}

As mentioned above, $S$-arithmetic groups over number fields are of type $F_{\infty}$. 
In contrast, $S$-arithmetic groups over function fields are almost never of type $F_{\infty}$. 

Conjecture 1(ii) implies that the upper bound in Conjecture 2 is sharp also in the number field case.
This indicates that are deep similarities from the point of view of 
 coarse geometry between $S$-arithmetic groups over function fields and those
over number fields:
In both cases the distortion dimension should be equal to  the Euclidean rank minus 1.

 A weakened and special case of Conjecture 2
has  recently been established in \cite{BEW}. Further evidence comes from the solution of Thurston's
conjecture, that $SL_n(\mathbb Z)$ has a quadratic $1$-dimensional Dehn function, by Young 
(see \cite{Yo4} for $n\geq 5$). Work of Lubotzky-Mozes-Raghunathan \cite{LMR} asserts that,
for $r\geq 2$,  $\mathbf G(\mathcal O_S)$ is undistorted in $\mathbf G_S$ in dimension $1$. Since this is not the case for $r=1$, it follows that 
$\textup{dist-dim}(\mathbf G(\mathcal O_S))=0$ if and only if the Euclidean rank of $\mathbf G_S$
is $1$. Moreover, results of Leuzinger-Pittet \cite{LP} (for number fields), Behr \cite{Be}
(for function fields) and Taback \cite{Ta} (for $PSL_2(\mathbb Z[\frac{1}{p}])$) show that
$\textup{dist-dim}(\mathbf G(\mathcal O_S))=1$, if the Euclidean rank of $\mathbf G_S$
is $2$.

Finally, we remark  that in view of Theorem 1 and Theorem \ref{euclower}, in order to prove Conjecture 2,  the problem is to  establish Euclidean isoperimetric inequalities below the 
rank.

\vspace{4ex}

{\bf Acknowledgement}. \ I thank Robert Young for pointing out to me some inaccuracies in an earlier version of this paper.

\vspace{4ex}

\noindent\textsc{Institute for Algebra und Geometry\\
 Karlsruhe Institute of Technology (KIT),  Germany}

\vspace{1ex}

\noindent{\tt enrico.leuzinger@kit.edu}
\end{document}